\documentclass{elsart}
\usepackage{amssymb,amsfonts}
\newcommand{\trace}{\mathop{\rm Tr}\nolimits}

\newcommand{\eig}{\mathop{\rm Eig}\nolimits}

\newcommand{\twomat}[4]{\left(\begin{array}{cc}#1&#2\\#3&#4\end{array}\right)}
\newcommand{\twovec}[2]{\left(\begin{array}{c}#1\\#2\end{array}\right)}

\newcommand{\schatten}[2]{\left|\left|\,{#2}\,\right|\right|_{#1}}
\newcommand{\geo}{\,\,\#\,\,}

\newcommand{\cS}{{\mathcal S}} % \cal not known

\DeclareRobustCommand\openone{\leavevmode\hbox{\small1\normalsize\kern-.33em1}}
\newcommand{\id}{\mathrm{\openone}}

\newcommand{\be}{\begin{equation}}
\newcommand{\ee}{\end{equation}}
\newcommand{\bea}{\begin{eqnarray}}
\newcommand{\eea}{\end{eqnarray}}
\newcommand{\beas}{\begin{eqnarray*}}
\newcommand{\eeas}{\end{eqnarray*}}

\newtheorem{theorem}{Theorem}
\newtheorem{lemma}{Lemma}
\newtheorem{corollary}{Corollary}

\newtheorem{proposition}{Proposition}
\newcount\minute
\newcount\hour
\def\currenttime{%
    \minute\time
    \hour\minute
    \divide\hour60
    \the\hour:\multiply\hour60\advance\minute-\hour\the\minute}
\begin{document}
\begin{frontmatter}
\title{A Norm Compression Inequality for Block Partitioned Positive Semidefinite Matrices}
\author{Koenraad M.R. Audenaert}
\address{
Imperial College London, The Blackett Lab--QOLS,
Prince Consort Road, London SW7 2BW, United Kingdom}
\address{
Institute for Mathematical Sciences, Imperial College London,
Exhibition Road, London SW7 2BW, United Kingdom}
\ead{k.audenaert@imperial.ac.uk}
\date{\today, \currenttime}
%------------------------------------------------------------------ ABSTRACT
\begin{abstract}
Let $A$ be a positive semidefinite matrix, block partitioned as
$$
A=\twomat{B}{C}{C^*}{D},
$$
where $B$ and $D$ are square blocks.
We prove the following inequalities for the Schatten $q$-norm $||.||_q$, which are sharp when the blocks are of
size at least $2\times2$:
$$
||A||_q^q \le (2^q-2) ||C||_q^q + ||B||_q^q+||D||_q^q, \quad 1\le q\le 2,
$$
and
$$
||A||_q^q \ge (2^q-2) ||C||_q^q + ||B||_q^q+||D||_q^q, \quad 2\le q.
$$
These bounds can be extended to symmetric partitionings into larger numbers of blocks, at the expense of no longer being sharp:
$$
||A||_q^q \le \sum_{i} ||A_{ii}||_q^q + (2^q-2) \sum_{i<j} ||A_{ij}||_q^q, \quad 1\le q\le 2,
$$
and
$$
||A||_q^q \ge \sum_{i} ||A_{ii}||_q^q + (2^q-2) \sum_{i<j} ||A_{ij}||_q^q, \quad 2\le q.
$$
\end{abstract}

\end{frontmatter}
%%%%%%%%%%%%%%%%%%%%%%%%%%%%%%%%%%%%%%%%%%%%%%%%%%%%%%%%%%%%%%%%%%%%%%%%%%%%%%%%%%%%%%%%%%%%%
\section{Introduction\label{sec1}}

%[TODO]
%Add further research directions:
%- bounds on $q$-norm of matrix given $q'$ norms of blocks, $q'\neq q$.
%- sharp bound for larger block numbers
%- Extensions of King to larger block numbers

In \cite{bhatia_kittaneh}, Bhatia and Kittaneh proved a number of interesting inequalities relating the
Schatten norms of a block partitioned operator to the Schatten norms of its constituent blocks.
Let the operator $T$ be written in block-matrix form as $T=[T_{ij}]$, with $1\le i,j\le d$,
then it is proven that, for example,
\be\label{b1}
d^{2-q} ||T||_q^q \le \sum_{i,j} ||T_{ij}||_q^q \le ||T||_q^q,\quad 2\le q
\ee
and
\be\label{b2}
d^{2-q} ||T||_q^q \ge \sum_{i,j} ||T_{ij}||_q^q \ge ||T||_q^q,\quad 1\le q\le 2.
\ee
It is also shown there that these inequalities are sharp.

In the following a bound will be called \textit{sharp} when it can be saturated for any allowed choice
of the constituent quantities of the bound. For example, in (\ref{b1}), these quantities are the norms
of the blocks $||T_{ij}||_q$, and sharpness means here that for any set of non-negative scalars $t_{ij}$
an operator $T$ exists such that $||T_{ij}||_q = t_{ij}$ and $||T||_q^q = \sum_{i,j}t_{ij}^q$.
Phrased differently, a sharp bound is the best possible bound exploiting \textit{a priori} specified knowledge.
This notion of sharpness is stronger than the one used in \cite{bhatia_kittaneh}.
Nevertheless, the second inequality in both (\ref{b1}) and (\ref{b2}) is evidently sharp according to our definition as well,
as can be seen by taking a $T$ with blocks $T_{ij}=[t_{ij}]\oplus 0$.

Inequalities like (\ref{b1}) and (\ref{b2}) are sometimes called \textit{norm compression inequalities},
because the full information contained in the operator is compressed into a smaller set of quantities,
the norms of its blocks, and the inequalities give useful bounds on the norm of the full operator
when only its compression is known.

In the present work we restrict attention to positive semidefinite (PSD) matrices. Under this extra restriction
bounds (\ref{b1}) and (\ref{b2}) are no longer sharp.
Indeed, by just considering the case $q=1$, which for positive matrices yields nothing but the trace,
we know that $||T||_1 = \sum_i ||T_{ii}||_1$, and the off-diagonal blocks should not contribute at all.

Known bounds of this form for PSD matrices and operators can be found in \cite{bhatia,HJII} and \cite{king}.
The best-known norm compression inequality (although it does not directly appear as such)
is probably the pinching inequality \cite{bhatia},
which holds for any weakly unitarily invariant norm, and arbitrary self-adjoint operators:
for any block-partitioned self-adjoint operator $A=[A_{ij}]\ge0$,
    \be
    |||A||| \ge |||\oplus_{i=1}^d A_{ii}|||.
    \ee
For Schatten norms, this reduces to
\be\label{pinch}
||A||_q \ge \left(\sum_{i=1}^d ||A_{ii}||_q^q\right)^{1/q},
\ee
which is indeed a norm compression inequality.
In (\cite{HJII}, p.~217 Problem 22)
one can find a complementary inequality for PSD $2\times 2$ block matrices, also valid for any unitarily invariant norm,
and readily extendible to PSD $d\times d$ block matrices:
\be\label{eq:hjii}
|||A||| \le \sum_{i=1}^d |||A_{ii}|||.
\ee
Here $|||A_{ii}|||$ is actually a shorthand for $|||A_{ii}\oplus 0|||$. That is,
the blocks have been implicitly filled out with zeroes to the same size as $A$.
There is a very simple proof of this inequality that also extends to operators:

    \textit{Proof.}
    Consider the $d=2$ case only. The general case follows by repartitioning the blocks iteratively.
    Fixing the diagonal blocks $A_{11}$ and $A_{22}$ fixes the RHS of (\ref{eq:hjii}),
    and restricts $A$ to a convex set whose extremal points are of the form $aa^*$, with $a=\twovec{a_1}{a_2}$
    and $a_i a_i^*=A_{ii}$. Here $a_1$ and $a_2$ are blocks with an equal number of columns.
    Because a norm, just as any convex function, reaches its maximum over a convex set in an extremal point of that set,
    we only need to check (\ref{eq:hjii}) for the extremal $A=aa^*$.
    Using the triangle inequality for norms, and the fact that $aa^*$ is unitarily equivalent with $a^*a\oplus 0$, we indeed get:
    \beas
    |||A||| &=& |||aa^*||| = |||a^* a||| = |||\sum_{i=1}^2 a_i^* a_i ||| \\
            &\le& \sum_{i=1}^2 |||a_i^* a_i||| = \sum_{i=1}^2 |||a_i a_i^*||| = \sum_{i=1}^2 |||A_{ii}|||.
    \eeas
    \qed

Bounds (\ref{pinch}) and (\ref{eq:hjii}) are sharp when the $q$-norms of the diagonal blocks only are known.
They are no longer sharp when the $q$-norms of all blocks are known, as can be seen by considering the
Frobenius norm (Schatten 2-norm).
Indeed, for that norm all blocks contribute evenly, while (\ref{pinch}) and (\ref{eq:hjii}) only take the diagonal blocks into
account.

What we are looking for in this paper are sharp norm compression inequalities for the Schatten norms of PSD block matrices,
when the norms of all the blocks are known, and not just the diagonal blocks.
Bounds of this kind have been discovered and proven by King \cite{king} for PSD $2\times2$ block matrices:
\be\label{king1}
\schatten{q}{\twomat{A_{11}}{A_{12}}{A_{21}}{A_{22}}}
\ge \schatten{q}{\twomat{||A_{11}||_q}{||A_{12}||_q}{||A_{21}||_q}{||A_{22}||_q}}, \quad 1\le q\le 2,
\ee
and
\be\label{king2}
\schatten{q}{\twomat{A_{11}}{A_{12}}{A_{21}}{A_{22}}}
\le \schatten{q}{\twomat{||A_{11}||_q}{||A_{12}||_q}{||A_{21}||_q}{||A_{22}||_q}}, \quad 2\le q.
\ee
That these bounds are sharp is easily seen by considering blocks $A_{ij}$ of the form $A_{ij}=[a_{ij}]\oplus 0$,
where $a_{ij}$ are non-negative scalars such that $a_{12}=a_{21}$ and $a_{11}a_{22}\ge a_{12}^2$.
In fact, when the $A_{ij}$ are scalars, equality holds in (\ref{king1}) and (\ref{king2}) throughout.

The obvious generalisation of (\ref{king1}) and (\ref{king2}) to higher numbers of blocks does not hold
for arbitrary $q$, although King has shown that $||A||_q\le\schatten{q}{(||A_{ij}||_q)_{i,j}}$
holds for \textit{integer} $q$ and any partitioning \cite{king2}.
For non-integer $q$ there are already counterexamples when the blocks $A_{ij}$ are scalars,
in which case the norm-compression is just the elementwise absolute value, which we denote here by $|A|$.
For example, for the matrix
$$
A=\left(\begin{array}{rrrr}
     2  &   0 &   -2  &  -2 \\
     0  &   2 &    2  &  -1 \\
    -2  &   2 &    3  &   0 \\
    -2  &  -1 &    0  &   2
\end{array}
\right)
$$
one finds $\schatten{1.5}{A}=7.7617$ and $\schatten{1.5}{|A|}=7.9761$.
We have not been able to find counterexamples for $3\times3$ partitionings, so it might be that (\ref{king1}) and (\ref{king2})
still hold in that case.

The underlying reason for the failure of (\ref{king1}) and (\ref{king2}) in the general case
seems to be that a norm compression maps a matrix to an elementwise non-negative
matrix. The natural ordering for those matrices is the elementwise ordering rather than the PSD ordering.
Likewise, unitarily invariant norms, which involve the eigenvalues of the matrix, do not seem to be the most natural
choice for norm compressions.
That King's bounds can be formulated for $2\times2$ (and may be $3\times3$) partitionings using
unitarily invariant norms is most likely a coincidence.

The main result of the present paper is a set of sharp bounds that is complementary to (\ref{king1}) and (\ref{king2}).
That is, for $1\le q\le 2$ we find an upper bound, and for $q\ge 2$ a lower bound on the $q$-norm of a $2\times2$ partitioned
PSD matrix, given the $q$-norms of its blocks.
These bounds are presented in Section \ref{sec3}.
In contrast to the bounds (\ref{king1}) and (\ref{king2}), our bounds can easily be generalised to any symmetric partitioning,
albeit at the expense of loss of sharpness.

Norm compression inequalities feature in proofs of the multiplicativity property
of the $1\rightarrow q$ norm of certain classes of completely positive maps.
Letting $\Phi$ be a completely positive (CP) map, this norm is defined as \cite{AH}
\be\label{norm}
||\Phi||_{1\rightarrow q} = \max_{||X||_1=1} ||\Phi(X)||_q,
\ee
where $X$ is Hermitian.
Multiplicativity of this norm w.r.t.\ the tensor product is the
statement that, for two CP maps $\Phi_1$ and $\Phi_2$ \cite{AH,amosov00}:
\be\label{mult}
||\Phi_1\otimes\Phi_2||_{1\rightarrow q} = ||\Phi_1||_{1\rightarrow q}\,\,||\Phi_2||_{1\rightarrow q}.
\ee
This basically says that the maximum in (\ref{norm}) for $\Phi=\Phi_1\otimes\Phi_2$ is achieved
for $X=X_1\otimes X_2$, where $X_i$ achieves the maximum in (\ref{norm}) for $\Phi_i$.
Multiplicativity (\ref{mult}) has been shown for various special classes of CP maps within
various ranges of $q$.
Unfortunately, there exists a class of channels for which (\ref{mult}) does not hold when $q>4.79$ \cite{werner02}.
Despite this counterexample to the general statement, (\ref{mult}) might still be true for any
tensor product of CP maps for values of $q$ close to 1.
If this were true, one could prove additivity of an entropic counterpart of (\ref{mult}), and with it a host
of other additivity results concerning CP maps. That would solve a number of long-standing open
problems in quantum information theory \cite{ka,shor}.
We intend to investigate the usefulness of our results in that setting in future work.
%%%%%%%%%%%%%%%%%%%%%%%%%%%%%%%%%%%%%%%%%%%%%%%%%%%%%%%%%%%%%%%%%%%%%%%%%%%%%%%%%%%%%%%%%%%%%%%%%%%%%%%%%%
\section{Preliminaries\label{sec2}}
%Schatten norms
The Schatten $q$-norms, for $1\le q<\infty$, are the non-commutative generalisation of the $l_q$ norms.
For a general matrix or operator $A$,
$$
||A||_q = (\trace(|A|^q))^{1/q},
$$
which reduces for positive semidefinite matrices $A$ to
$$
||A||_q = (\trace(A^q))^{1/q}.
$$

We will use the positive semidefinite ordering on Hermitian matrices throughout, denoted
$A\ge B$, which means that $A-B\ge 0$.
This ordering is preserved under arbitrary conjugations: $A\ge B$ implies $XAX^*\ge XBX^*$ for arbitrary $X$.

%PSD condition 2X2
It is well-known that a $2\times 2$ block-matrix $A=\twomat{B}{C}{C^*}{D}$
with positive definite $B$ and $D$ is positive semidefinite if and only if $B\ge CD^{-1}C^*$.

% geometric mean
The set $S$ of Hermitian $C$ such that $\twomat{B}{C}{C}{D}$ is PSD, has a unique maximum, called
the \textit{geometric mean} of $B$ and $D$ \cite{kubo_ando,pusz_woronowicz}.
For any $A,B>0$, the geometric mean of $A$ and $B$, denoted $A\geo B$, is given by
\be\label{geodef}
A\geo B = B\geo A = A^{1/2}(A^{-1/2}BA^{-1/2})^{1/2} A^{1/2}.
\ee
For $A,B\ge 0$, the geometric mean is defined by
$$
A\geo B = \lim_{\epsilon\downarrow 0}(A+\epsilon\id)\geo(B+\epsilon\id).
$$
For $A$ and $B$ commuting, (\ref{geodef}) reduces to $A\geo B = (AB)^{1/2}$.

As basic properties, we need \cite{ando,ando_hiai}:
\begin{itemize}
\item $C(A\geo B)C^* = (CAC^*)\geo(CBC^*)$;
\item $(A\geo B)^{-1} = A^{-1} \geo B^{-1}$;
\item $(A,B)\mapsto A\geo B$ is jointly monotone in its arguments.
That is: if $A_1\le A_2$ and $B_1\le B_2$, then also $A_1\geo B_1 \le A_2\geo B_2$.
\end{itemize}
We will also need the following Lemma:
\begin{lemma}\label{lem1}
For $A,B>0$, the unique positive definite solution of the equation
$XA^{-1}X=B$ is given by $X=A\geo B$.
\end{lemma}
\textit{Proof.}
From $XA^{-1}X=B$ it follows that $X$ is in the set $S$ of Hermitian matrices $C$ for which $\twomat{A}{C}{C}{B}\ge0$,
hence $X\le A\geo B$.
It also follows that $X^{-1}AX^{-1}=B^{-1}$, hence $X^{-1}\le A^{-1}\geo B^{-1} = (A\geo B)^{-1}$.
Thus, if we restrict to positive definite $X$, we find $X\ge A\geo B$.
Therefore, we actually have equality: $X=A\geo B$.
\qed

%power means
A generalisation of the geometric mean is the $\alpha$-\textit{power mean}, for $0\le\alpha\le 1$ and $A,B>0$:
$$
A\,\,\#_\alpha\,\, B = A^{1/2}(A^{-1/2}BA^{-1/2})^\alpha A^{1/2}.
$$

%operator mono,conc,conv
%operations that preserve that
A matrix function $f$ is \textit{operator monotone} iff it preserves the PSD ordering, i.e.\ $A\ge B$ implies $f(A)\ge f(B)$.
If $A\ge B$ implies $f(A)\le f(B)$, we say $f$ is \textit{inversely operator monotone}.
A matrix function $f$ is \textit{operator convex} iff for all $0\le \lambda\le 1$ and for all $A,B\ge 0$,
$$
f(\lambda A+(1-\lambda)B)\le \lambda f(A)+(1-\lambda)f(B).
$$
If $-f$ is operator convex, we say $f$ is \textit{operator concave}.

The primary matrix function $x\mapsto x^p$ is
operator convex for $1\le p\le 2$,
operator monotone and operator concave for $0\le p\le 1$,
and inversely operator monotone and operator convex for $-1\le p\le 0$ \cite{bhatia}.

%log-majo
We will also make use of the log-majorisation relation for positive $A$, $B$:
$$
A \prec_{\log} B \Longleftrightarrow \log A\prec \log B,
$$
which implies weak majorisation $A\prec_w B$, and hence $|||A|||\le |||B|||$ for any unitarily
invariant norm.

%metric
Finally, we will use the $\delta_\infty$ metric on the positive cone, defined as
$$
\delta_\infty(A,B) = ||\log\eig(AB^{-1})||_\infty,
$$
for $A,B>0$.
Here, $\eig(A)$ is the vector of eigenvalues of $A$,
and the norm used is the $l_\infty$ vector norm.
This metric is well-defined since, for $A,B>0$, $AB^{-1}$ has positive eigenvalues.
We note that
$$
\delta_\infty(A,B) = \max(|\log\lambda_1^\downarrow(AB^{-1})|, |\log\lambda_1^\uparrow(AB^{-1})| ),
$$
where $\lambda_1^\downarrow$ and $\lambda_1^\uparrow$ denote the largest and smallest eigenvalue, respectively.
%%%%%%%%%%%%%%%%%%%%%%%%%%%%%%%%%%%%%%%%%%%%%%%%%%%%%%%%%%%%%%%%%%%%%%%%%%%%%%%%%%%%%%%%%%%%%%%%%%%%%%%%%%
\section{Main Result\label{sec3}}
\begin{theorem}
Let $A$ be a positive semidefinite block matrix
$$
A=\twomat{B}{C}{C^*}{D},
$$
where $B$ and $D$ are square blocks.
Then we have the following bound on the Schatten $q$-norm of $A$ for $1\le q\le 2$:
\be\label{ineq}
||A||_q^q \le (2^q-2) ||C||_q^q + ||B||_q^q+||D||_q^q.
\ee
\end{theorem}
It is easy to see that, for $q=1$ and for $q=2$, equality holds.
Indeed, for $q=1$, (\ref{ineq}) reduces to $\trace(A)=\trace(B)+\trace(C)$, and
for $q=2$, $\trace(A^2) = \trace(B^2)+2\trace(|C|^2)+\trace(D^2)$.
In this sense, (\ref{ineq}) interpolates between these two extremal cases.

Using a standard duality argument, we find that for $q\ge 2$, inequality (\ref{ineq}) is reversed:
\begin{corollary}\label{co:rev}
For $q\ge 2$, and with $A$, $B$, $C$, $D$ as in Theorem 1,
\be\label{ineqrev}
||A||_q^q \ge (2^q-2) ||C||_q^q + ||B||_q^q+||D||_q^q.
\ee
\end{corollary}
\textit{Proof.}
Consider the matrix $A=\twomat{B}{C}{C^*}{D}$ from Theorem 1.
We will restrict attention to the case where $B$ and $D$ are of equal size, so that $C$ is square.
Evidently, the blocks can always be filled out with zeroes to bring them to this form without changing
the validity of the bound.
Furthermore, we restrict to $C=C^*$ that are positive semidefinite.
To see that this incurs no loss of generality either,
consider the polar decomposition of general $C$, $C=UC'$, where $U$ is a unitary and $C'\ge0$.
Then
$$
A':=\twomat{U^*}{0}{0}{\id} \, A \, \twomat{U}{0}{0}{\id} = \twomat{B'}{C'}{C'}{D},
$$
with $B'=U^*BU$.
Clearly, $A$ and $A'$ have the same norm, and so do $B$ and $B'$, and $C$ and $C'$.
Therefore, in the following, we can take $C\ge0$, so that all occurrences of $||.||_q^q$ can be written as $\trace(.)^q$.

Let $q\ge2$ and let $p$ be the conjugate power of $q$: $1/p+1/q=1$.
H\"older's inequality for positive semidefinite $A$ and $B$ reads
$\trace[AB]\le ||A||_p\,||B||_q$, with equality if $B=A^{p-1}$.
This allows one to express the norm $||A||_p$ as the supremum of $\trace[AB]$ over all $B\ge0$ for which $||B||_q=1$.
In other words, for every $A\ge0$ there exists an optimal $B\ge0$ with $||B||_q=1$ such that $||A||_p=\trace[AB]$,
and for all other $B\ge0$ with $||B||_q=1$ one has $||A||_p\ge\trace[AB]$.
As the optimal $B$ is given by $A^{p-1}/\schatten{q}{A^{p-1}}$, one can always safely assume
that the optimal $B$ has the same direct sum structure as $A$ has.

Now consider the expression
\be\label{bdc}
||B\oplus D\oplus (2^p-2)^{1/p} C||_p.
\ee
Let $P$, $Q$ and $R$ be positive semidefinite matrices such that
$P\oplus R\oplus (2^q-2)^{1/q} C$ is optimal for the norm in (\ref{bdc}) in the abovementioned sense.
That is:
$$
||P\oplus R\oplus (2^q-2)^{1/q} C||_q = 1,
$$
and
\beas
\lefteqn{||B\oplus D\oplus (2^p-2)^{1/p} C||_p} \\
&=& \trace\left[(B\oplus D\oplus (2^p-2)^{1/p} C)\,(P\oplus R\oplus (2^q-2)^{1/q} C)\right] \\
&=& \trace[BP+DR+(2^p-2)^{1/p}(2^q-2)^{1/q}CQ].
\eeas
Now notice that for all $q$, $(2^p-2)^{1/p}(2^q-2)^{1/q}\le 2$, with equality in $q=2$.
Thus
\beas
\schatten{p}{B\oplus D\oplus (2^p-2)^{1/p} C} &\le& \trace[BP+DR+2CQ] \\
&=& \trace\left[\twomat{B}{C}{C}{D}\,\twomat{P}{Q}{Q}{R}\right].
\eeas
On the other hand, from $\schatten{q}{P\oplus R\oplus (2^q-2)^{1/q} C} = 1$ and Theorem 1, it follows that
$\schatten{q}{\twomat{P}{Q}{Q}{R}}\le 1$.
Thus, using H\"older's inequality,
we may conclude that $\trace\left[\twomat{B}{C}{C}{D}\,\twomat{P}{Q}{Q}{R}\right] \le \schatten{p}{\twomat{B}{C}{C}{D}}$,
which proves the inequality (\ref{ineqrev}) of the Corollary.
\qed

We can combine (\ref{ineq}) with (\ref{b2}), applied to the $C$ block,
to generalise our bounds to general $d\times d$ partitionings, by repartitioning the $B$ and $C$ blocks recursively.
\begin{corollary}
For any PSD matrix $A$, partitioned into $d\times d$ blocks $A_{ij}$ such that the diagonal blocks are square,
\be\label{gen}
||A||_q^q \le \sum_i ||A_{ii}||_q^q + (2^q-2) \sum_{i<j}||A_{ij}||_q^q, \quad 1\le q\le 2
\ee
and
\be\label{genrev}
||A||_q^q \ge \sum_i ||A_{ii}||_q^q + (2^q-2) \sum_{i<j}||A_{ij}||_q^q, \quad 2\le q.
\ee
\end{corollary}
The proof of (\ref{ineqrev}) extends without essential changes to (\ref{genrev}).

Concerning sharpness, we first have to mention that for blocks of size $1\times 1$, our bounds are not sharp, quite simply because
King's bounds (\ref{king1}) and (\ref{king2}) are equalities in that case. For blocks of size $2\times 2$ (and larger),
our bounds (\ref{ineq}) and (\ref{ineqrev}) are sharp, as witnessed by blocks of the form
$$
B = \twomat{b}{0}{0}{c},\quad
C = \twomat{0}{0}{0}{c},\quad
D = \twomat{d}{0}{0}{c},
$$
where $b$, $c$ and $d$ are non-negative numbers.
In Section \ref{sec:notsharp}, however, we show that (\ref{gen}) is not sharp.
It would be interesting to find better bounds for that case, but at this point it is not clear to us
whether this question has a reasonable answer.

To prove the central technical result (\ref{ineq}), we can, just as in the proof of Corollary \ref{co:rev},
w.l.o.g.\ restrict attention to the case where block $C$ is square and positive semidefinite.
Inequality (\ref{ineq}) can then be reformulated in a way that sheds light on the somewhat curious factor of $2^q-2$.
Note, namely, that
$$
\trace\twomat{C}{C}{C}{C}^q = 2^q\trace C^q,
$$
and
$$
\trace\twomat{C}{0}{0}{C}^q = 2\trace C^q.
$$
Hence, (\ref{ineq}) can be written as
\be\label{ineq2}
\trace\twomat{B}{C}{C}{D}^q - \trace\twomat{B}{0}{0}{D}^q \le \trace\twomat{C}{C}{C}{C}^q - \trace\twomat{C}{0}{0}{C}^q.
\ee
It is clear that both sides are non-negative, since $\twomat{B}{0}{0}{D}$ is a pinching of $\twomat{B}{C}{C}{D}$,
and weakly unitarily invariant norms, such as the Schatten norms, are non-increasing under pinchings.
The difference expressed by the left-hand side is thus the amount of norm decrease caused by this particular pinching,
and the inequality says that, when fixing $C$ and constraining $B$ and $D$ to keep $A$ PSD,
this norm decrease is maximal when $B=D=C$.

The Proof of Theorem 1 will be given in Sections \ref{sec4} to \ref{sec7}.
%%%%%%%%%%%%%%%%%%%%%%%%%%%%%%%%%%%%%%%%%%%%%%%%%%%%%%%%%%%%
\section{Bound (\ref{gen}) is not sharp}\label{sec:notsharp}
In this Section, we consider the generalisation (\ref{gen}) of our bound to general partitionings, and show that
it is no longer sharp.
We consider a particular class of PSD matrices $(A_{i,j})_{i,j}$ for which every block has the same $q$-norm:
$||A_{ij}||_q=a$. We first show that this implies that all blocks have the same absolute value.

Consider the blocks $A_{ii}$, $A_{jj}$ and $A_{ij}$ for some $i<j$.
Non-negativity of $A$ implies that $A_{ii}\ge A_{ij}A_{jj}^{-1}A_{ij}^*$.
Since all blocks have the same norm, we actually must have equality.
\begin{lemma}
For a PSD block matrix $A=\twomat{B}{C}{C^*}{D}\ge 0$, the equality $||B||_q=||C||_q=||D||_q$
implies $B = CD^{-1}C^*$.
\end{lemma}
\textit{Proof.}
Suppose there was a $\Delta\ge 0$ for which $B = CD^{-1}C^*+\Delta$.
By \cite{bhatia}, (IV.53), $|||(A+B)\oplus 0|||\ge |||A\oplus B|||$ for $A,B\ge 0$, hence
$\trace(A+B)^q \ge \trace A^q+\trace B^q$. For finite $q$ this means that $||A+B||_q$ is strictly larger than $||A||_q$ when
$B$ is non-zero. Specifically, if $\Delta$ is non-zero, we find $||B||_q > ||CD^{-1}C^*||_q$.
Using a Theorem of Horn and Mathias \cite{hm}, $||CD^{-1}C^*||_q\ge ||C||_q^2 ||D||_q^{-1}$,
hence the non-vanishing of $\Delta$ implies $||CD^{-1}C^*||_q > ||C||_q^2 ||D||_q^{-1}$, which violates the
statement that $||B||_q=||C||_q=||D||_q=a$. Therefore, $\Delta$ must be zero.
\qed

Using King's inequality (\ref{king1}), we can strenghten this further.
\begin{lemma}
For a PSD block matrix $A=\twomat{B}{C}{C^*}{D}\ge 0$, the equality $||B||_q=||C||_q=||D||_q=a$, $1< q\le 2$,
implies $B = D = UC$, where $U$ is a unitary commuting with $D$.
Thus, in some basis, $B$, $C$ and $D$ are diagonal, and $B=|C|=D$.
\end{lemma}
\textit{Proof.}
From the previous Lemma, we already know that $B = CD^{-1}C^*$.
Using (\ref{king1}), we find
$$
\schatten{q}{\twomat{CD^{-1}C^*}{C}{C^*}{D}} \ge \schatten{q}{\twomat{a}{a}{a}{a}} = 2a.
$$
The left-hand side is equal to $||D+D^{-1/2}C^*CD^{-1/2}||_q$. By the triangle inequality,
\beas
||D+D^{-1/2}C^*CD^{-1/2}||_q &\le& ||D||_q + ||D^{-1/2}C^*CD^{-1/2}||_q \\
&=& ||D||_q + ||CD^{-1}C^*||_q = ||D||_q+||B||_q = 2a.
\eeas
Combining these two inequalities, we find that equality holds. Now, by the Lemma below, this implies
$D=tB$, with, in particular, $t=1$, thus $D=B$.
This further implies $CD^{-1}C^*=D$ and also $D^{-1/2}C^*CD^{-1/2}=D$.
From the latter equation we find $|C|=D$. The polar decomposition of $C$ must therefore be $C=UD$.
Inserting this in the former equation yields $UDU^*=D$, so that $U$ must commute with $D$.
\qed

\begin{lemma}
For given matrices $A,B$,
equality in the Triangle Inequality for $q$-Schatten norms with $1<q\le 2$,
$$
||A+B||_q = ||A||_q+||B||_q,
$$
implies $A=tB$, for some $t\ge0$.
\end{lemma}
\textit{Proof.}
By convexity of norms, for all $\lambda\in[0,1]$,
$$
||\lambda A+(1-\lambda)B||_q \le \lambda ||A||_q+(1-\lambda)||B||_q.
$$
Then $||A+B||_q = ||A||_q+||B||_q$ implies equality for all $\lambda$, and by dividing both sides by $\lambda$, we get
$$
|| A+tB||_q \le ||A||_q+t||B||_q,
$$
where $t=(1-\lambda)/\lambda>0$.
Choosing $t$ equal to $||A||_q/||B||_q$ and setting $B'=tB$, we get, in particular,
$||A||_q=:a$, $||B'||_q=a$, and $||A+B'||_q = 2a$.
Inserting this in the ``hard'' Clarkson-McCarthy inequality \cite{simon}, which is valid for $1\le q\le 2$:
$$
||A+B'||_q^p + ||A-B'||_q^p \le 2(||A|_q^q+||B'||_q^q)^{p/q},
$$
with $1/p+1/q=1$, gives, for $q>1$ (i.e.\ finite $p$)
$$
||A-B'||_q^p \le (2\,2^{p/q}-2^p) a^p = 0,
$$
whence it follows that $A=tB$.
\qed

So, we now can already conclude that $A=(A_{ij})_{ij}$ must in a certain basis be of the form
$$
A=\bigoplus_j x_j X_j,
$$
with $x_j\ge0$ such that $\sum_j x_j^q=a^q$,
and $X_j$ $d\times d$ PSD matrices whose elements all have modulus 1.
Now, this can only be if the $X_j$ are rank 1,
as can be seen by noting that $\trace (X_j/d)^2=\trace(X_j/d)$.
Thus, $||A||_q^q = \sum_j x_j^q d^q = a^q d^q$.
On the other hand, (\ref{gen}) gives $||A||_q^q \le (2^q-2) (d(d-1)/2)a^q + d a^q$.
As this is strictly larger than $a^q d^q$ for $1<q<2$, this shows that (\ref{gen}) is not sharp.
%%%%%%%%%%%%%%%%%%%%%%%%%%%%%%%%%%%%%%%%%%%%%%%%%%%%%%%%%%%%%%%%%%%%%%%%%%
\section{Proof of Theorem 1\label{sec4}}
We only have to prove (\ref{ineq2}) for $1<q<2$. The cases $q=1$ and $q=2$ are trivial, as noted before.
Furthermore, we only have to deal with the case where all blocks are square and of the same size,
We can easily generalise our Main Theorem to non-square $C$ blocks,
by filling out the smaller blocks with zeroes to the required size.

We deal first with the case that $B$ and $D$ are bounded and positive definite,
and leave the remaining cases for last (cfr. Proposition \ref{prop2}).

Let us consider the left-hand side of (\ref{ineq2}) and effectively calculate its maximum value.
We start by maximising it over $B$.
The constraint on $B$, originating from the requirement $A\ge0$, is $B\ge CD^{-1}C$.
We will now show that the maximum over $B$ is obtained in $B=B_0:=CD^{-1}C$.
Let us thereto put $B=B_0+t\Delta$, with $\Delta\ge0$, and define
$$
f(t) := \trace\twomat{B_0+t\Delta}{C}{C}{D}^q - \trace\twomat{B_0+t\Delta}{0}{0}{D}^q.
$$
The derivative of $f$ is given by
$$
f'(t) = q\trace\left[\left(\twomat{B}{C}{C}{D}^{q-1}-\twomat{B}{0}{0}{0}^{q-1}\right)\,\twomat{\Delta}{0}{0}{0}\right].
$$
Introducing the projector $P=\id\oplus 0$, we can write
$$
f'(t) = q\trace\left[\left(P\twomat{B}{C}{C}{D}^{q-1}P-\left(P\twomat{B}{C}{C}{D}P\right)^{q-1}\right)\,\twomat{\Delta}{0}{0}{0}\right].
$$
For $1< q\le 2$, the function $x\mapsto g(x)=x^{q-1}$ is operator concave on $[0,+\infty)$, and $g(0)=0$.
Therefore (\cite{bhatia}, Theorem V.2.3)
$$
P\twomat{B}{C}{C}{D}^{q-1}P \le \left(P\twomat{B}{C}{C}{D}P\right)^{q-1}.
$$
This shows that $f'(t)\le0$ and that $f(t)$ is indeed maximal in 0.
Therefore, we can henceforth put $B=CD^{-1}C$.

\bigskip
\bigskip
Define $f(D)$ as
\be\label{fd}
f(D) := \trace\twomat{CD^{-1}C}{C}{C}{D}^q - \trace\twomat{CD^{-1}C}{0}{0}{D}^q.
\ee
Since
\be\label{cd}
\twomat{CD^{-1}C}{C}{C}{D} = \twovec{CD^{-1/2}}{D^{1/2}}(D^{-1/2}C \,\,\, D^{1/2}),
\ee
and $CD^{-1}C$ has the same spectrum as $D^{-1/2}C^2D^{-1/2}$,
we can rewrite $f(D)$ as
\be
f(D) = \trace(G+D)^q - \trace G^q - \trace D^q,
\ee
where we have introduced
\be
G = D^{-1/2} C^2 D^{-1/2}.
\ee
A short (numerical) calculation reveals that $f(D)$ is neither convex nor concave, not even in the scalar case ($C$ and $D$ scalars).

To perform the maximisation of $f(D)$ over all possible $D>0$,
we calculate the gradient and stationary points of $f(D)$. We replace $D$ by $D+tX$, with Hermitian $X$,
and calculate the Fr\'echet derivative of (\ref{fd}):
\beas
\frac{\partial}{\partial t}\Big|_{t=0} f(D+tX) &=& q\,\trace\Big[
X\,D^{-1/2} \Big(
D((D+G)^{q-2}-D^{q-2})D \\
&& \qquad\qquad\qquad - G((D+G)^{q-2}-G^{q-2})G
\Big) D^{-1/2}
\Big].
\eeas
In this calculation we have used the approximation
\beas
(D+tX)^{-1} &=& D^{-1/2}(\id+t D^{-1/2}XD^{-1/2})^{-1}D^{-1/2} \\
&=& D^{-1/2}(\id-t D^{-1/2}XD^{-1/2})D^{-1/2} + O(t^2) \\
&=& D^{-1} - t D^{-1}XD^{-1} + O(t^2),
\eeas
the expression for the Fr\'echet derivative of the power function
$$
\frac{\partial}{\partial t}\Big|_{t=0} \trace(A+t\Delta)^q = q\trace(A^{q-1}\Delta),
$$
and the equality
$$
\twomat{CD^{-1}C}{C}{C}{D}^p = \twovec{CD^{-1/2}}{D^{1/2}}\, (D^{-1/2} C^2 D^{-1/2}+D)^{p-1}\, \Bigg(D^{-1/2}C\,\,\,\,D^{1/2}\Bigg),
$$
for all $p$, which follows from (\ref{cd}).
Therefore, the gradient of $f(D)$ is given by the expression
\bea
\nabla f(D) &=& qD^{-1/2}\big[D((D+G)^{q-2}-D^{q-2})D  \nonumber\\
&& \qquad\quad - G((D+G)^{q-2}-G^{q-2})G\big] D^{-1/2},
\eea
and $D$ is a stationary point of $f(D)$ if and only if this gradient is zero.
This clearly shows that the gradient of $f$ is well-defined and continuous in the interior of the positive semidefinite cone $\cS$.
It is also clear that $D=C$, implying that also $G=C$, is a stationary point.

The global maximum of $f$ must either be a stationary point, a singular point, or a boundary point.
As the gradient of $f$ is well-defined in the interior of $\cS$, $f$ has no singular points.
In the following Sections we prove that $D=G=C$ is the only stationary point of $f$.
More precisely, in Sections \ref{sec5} and \ref{sec6} we will prove the following Proposition:
\begin{proposition}\label{prop1}
For $p$ in the range $-1<p< 1$, $p\neq 0$, and for $D>0$, the equation in $G$
$$
D((D+G)^p-D^p)D - G((D+G)^p-G^p)G = 0
$$
has one solution over the positive definite matrices, namely $G=D$.
\end{proposition}
Since we are dealing with values $1<q<2$, this Proposition applies with $p=q-2$.

Finally, we show in Section \ref{sec7} that the values of $f$ on the boundary of $\cS$ are not greater than $f(C)$.
This is proven in an inductive way, as follows:
\begin{proposition}\label{prop2}
Assuming $f(D)\le f(C)$ holds for all $C$ and $D$ of size $d'\times d'$,
$f(D)\le f(C)$ also holds for $d\times d$ matrices $D$ that are bounded and invertible on a $d'$-dimensional subspace
of the full $d$-dimensional space (with $d'<d$).
\end{proposition}
Using induction on the size $d$ of the blocks,
these two Propositions allow us to conclude that $D=C$, the ``only stationary point in town'',
is the global maximum of $f(D)$, so that $f(D)\le f(C)$ for all $D\ge 0$, which is what we needed to show.
This finishes the proof of Theorem 1.
%%%%%%%%%%%%%%%%%%%%%%%%%%%%%%%%%%%%%%%%%%%%%%
\section{Uniqueness of the stationary point\label{sec5}}
In this and the following Section, we present the proof of Proposition \ref{prop1}.
We consider the equation
\be\label{eq}
D((D+G)^p-D^p)D = G((D+G)^p-G^p)G
\ee
over $G>0$,
and we will show that $G=D$, implying $G=D=C$, is its only solution for values of $p$, $-1< p< 1$, $p\neq 0$.

We start with the case $0< p< 1$.
Applying Lemma \ref{lem1}, (\ref{eq}) is equivalent with
$$
G = \left(D((D+G)^p-D^p)D\right) \geo \left((D+G)^p-G^p\right)^{-1},
$$
and we define the map $\Phi_D$ that maps $G$ to the matrix expressed by the right-hand side of this equation:
\be\label{map}
G\mapsto \Phi_D(G)= \left(D((D+G)^p-D^p)D\right) \geo \left((D+G)^p-G^p\right)^{-1}.
\ee
For the case $-1< p< 0$, $(D+G)^p-D^p$ and $(D+G)^p-G^p$ are negative, and we now find
$$
G = \left(D(D^p-(D+G)^p)D\right) \geo \left(G^p-(D+G)^p\right)^{-1}.
$$
The sign changes, as compared to (\ref{map}), are necessary for the geometric mean to have positive definite arguments.
Therefore, in that case, we define $\Phi_D$ as
\be\label{mapneg}
G\mapsto \Phi_D(G)= \left(D(D^p-(D+G)^p)D\right) \geo \left(G^p-(D+G)^p\right)^{-1}.
\ee

To prove that (\ref{eq}) has only one solution, we will show that $\Phi_D$ has only one fixed point (namely $G=D$)
for $-1<p< 1$, $p\neq 0$.
The way we will do this is by showing that $\Phi_D$ is ``contractive w.r.t.\ the fixed point $D$''.
Endowing the cone of positive semidefinite matrices $\cS$ with the metric $\delta_\infty$,
contractivity of $\Phi_D$ w.r.t.\ $D$ means the inequality
\be\label{contract}
\delta_\infty(\Phi_D(G),D)\le\beta \delta_\infty(G,D),
\ee
where the ``Lipschitz constant'' $\beta$ is strictly less than 1.
This statement resembles the definition of contractivity of a map, which says that, for all $G$ and $G'$,
$\delta_\infty(\Phi(G),\Phi(G'))\le\beta \delta_\infty(G,G')$, with Lipschitz constant $\beta<1$.
By the contraction mapping principle, contractive maps have a unique fixed point in $\cS$.
Similarly, the weaker statement (\ref{contract}) is already enough to show that $D$ is the unique fixed point of $\Phi_D$.
Indeed, suppose there is another fixed point $D'$: $\Phi_D(D')=D'$. Taking $G=D'$ in (\ref{contract})
then yields $\delta_\infty(D',D)\le\beta \delta_\infty(D',D)$, which can only be true if $\delta_\infty(D',D)=0$, i.e.\ $D'=D$.
%%%%%%%%%%%%%%%%%%%%%%%%%%%%%%%%%%%%%%%%%%%%%%%%%%%%%%%%%%%%%%%%%%%%%%%%%%%
\section{Contractivity of the map $\Phi_D$\label{sec6}}
We will now prove that when $-1\le p\le 1$, (\ref{contract}) holds with $\beta = p/(2^{p+1}-2)$, which is strictly less than 1 for $-1<p$.
If the map $\Phi_D$ would have been operator monotone, this would have allowed us to straightforwardly reduce the problem to the scalar case.
However, the subexpression $((D+G)^p-G^p)^{-1}$ is not monotone in $G$. Nevertheless, monotonicity holds in the
following very restricted sense, and this will turn out to be just enough for our purposes.
%----------------------------------------------------------- Lemma 2
\begin{lemma}\label{lem2}
Let $A,B$ be positive semidefinite and $k$ a positive scalar.

For $0\le p\le 1$:\\
$A\le kB$ implies $(A+B)^p-A^p \ge (kB+B)^p-(kB)^p \ge 0$.

For $-1\le p\le 0$, the orderings are reversed:\\
$A\le kB$ implies $(A+B)^p-A^p \le (kB+B)^p-(kB)^p \le 0$.
\end{lemma}
As a side remark, we note that, for instance for $0\le p\le 1$,
$A\ge kB$ does \textit{not} imply $(A+B)^p-A^p \le (kB+B)^p-(kB)^p$.

\textit{Proof.}
We note first that $A+B$ can be written as the convex combination $\lambda(k+1)B+(1-\lambda)((k+1)/k) A$,
with $\lambda=1/(k+1)$.

By operator concavity of the function $x\mapsto x^p$, $0\le p\le 1$, we then have
$$
(A+B)^p \ge \lambda (k+1)^p B^p+(1-\lambda)((k+1)/k)^p A^p,
$$
so that
$$
(kB+B)^p - (A+B)^p \le \left(\frac{k}{k+1}\right)^{1-p} ((kB)^p-A^p).
$$
Since $x\mapsto x^p$, $0\le p\le 1$, is also operator monotone, $(kB)^p-A^p \ge 0$.
For $p\le 1$ and $k\ge 0$, the factor $(k/(k+1))^{1-p}$ is $\le 1$, so that
$(kB+B)^p - (A+B)^p \le (kB)^p-A^p$ follows, which is equivalent to the first inequality of the Lemma.

For the second case, $-1\le p\le 0$, we proceed in exactly the same way, but now
exploiting the operator convexity and inverse monotonicity of $x\mapsto x^p$ for $-1\le p\le 0$.
\qed
%------------------------------------------------------------- End of Lemma 2

Using Lemma \ref{lem2}, we can easily prove similar statements for $\Phi_D(G)$.
Define the function
\be
\phi(x) = \Phi_1(x) = \left(\frac{(1+x)^p-1}{(1+x)^p-x^p}\right)^{1/2}.
\ee
It is readily seen that $\phi(1/x) = 1/\phi(x)$.
%-------------------------------------------------------------------------- Lemma 3
\begin{lemma}\label{lem3}
Consider matrices $D,G> 0$, and a scalar $k> 0$.
For $-1\le p\le 1$,
\bea
G\le k D&\mbox{ implies }&\Phi_D(G) \le \phi(k) D, \label{lem31}\\
D\le k G&\mbox{ implies }&\Phi_D(G) \ge \phi(k)^{-1} D. \label{lem32}
\eea
\end{lemma}
\textit{Proof.}
We start with the case $0\le p\le 1$, for which the function $x\mapsto x^p$ is operator monotone (and concave).
Then $G\le k D$ implies
\beas
D((D+G)^p-D^p)D &\le& D((D+kD)^p-D^p)D \\
&=& ((1+k)^p-1) D^{p+2}.
\eeas
By Lemma \ref{lem2}, we also have
\beas
((D+G)^p-G^p)^{-1} &\le& ((D+kD)^p-(kD)^p)^{-1} \\
&=& ((1+k)^p-k^p)^{-1} D^{-p}.
\eeas
Joint monotonicity of the geometric mean then yields
\beas
\Phi_D(G) &\le& ((1+k)^p-1) D^{q+1} \geo ((1+k)^p-k^p)^{-1} D^{-p} \\
&=& \phi(k) D,
\eeas
which is (\ref{lem31}).

To prove (\ref{lem32}), $D\le k G$ similarly implies
$$
((D+G)^p-G^p)^{-1} \ge ((1+k)^p-1)^{-1} G^{-p}.
$$
Using Lemma \ref{lem2} again, we have
\beas
D((D+G)^p-D^p)D &\ge& D((kG+G)^p-(kG)^p)D \\
&=& ((1+k)^p-k^p)DG^pD.
\eeas
For the geometric mean we get
\beas
\Phi_D(G) &\ge& ((1+k)^p-1)^{-1} G^{-p} \geo ((1+k)^p-k^p)DG^pD \\
&=& \phi(k)^{-1} (G^{-p} \geo DG^pD) \\
&=& \phi(k)^{-1} D.
\eeas
In the last line we have used
\beas
G^{-p} \geo DG^pD &=& D^{1/2}(D^{-1/2}G^{-p}D^{-1/2} \geo D^{1/2}G^pD^{1/2})D^{1/2} \\
&=& D^{1/2}((D^{1/2}G^pD^{1/2})^{-1} \geo D^{1/2}G^pD^{1/2})D^{1/2} \\
&=& D.
\eeas

For $-1\le p\le 0$, inequalities (\ref{lem31}) and (\ref{lem32}) are proven in exactly the same way. On one hand, since
$x\mapsto x^p$ is now inversely operator monotone, the inequality signs are reversed, and the same applies for
the inequality of Lemma \ref{lem2}.
However, this reversal is counteracted by the fact
that in this regime $\Phi_D(G)$ is defined by (\ref{mapneg}), which has additional sign changes,
hence the inequalities of the Lemma still remain valid.
\qed
%--------------------------------------------------------------------- End of Lemma 3

From this Lemma we get inequalities for $\lambda_1^\downarrow$ and $\lambda_1^\uparrow$
of $GD^{-1}$ and $\Phi_D(G)D^{-1}$, valid for $-1\le p\le 1$.
Assume first that $\lambda_1^\downarrow(G D^{-1})=K$. This amounts to $G\le KD$, and by the first statement of Lemma \ref{lem3},
implies $\Phi_D(G)\le \phi(K)D$, hence $\lambda_1^\downarrow(\Phi_D(G)D^{-1})\le \phi(K)$.
Thus we get
\be\label{eqlam1}
\lambda_1^\downarrow(\Phi_D(G)D^{-1})\le \phi(\lambda_1^\downarrow(G D^{-1})).
\ee
Then assume $\lambda_1^\uparrow(G D^{-1})=k$, which means
that $G\ge kD$, and by the second statement of Lemma \ref{lem3},
$\Phi_D(G)\ge (1/\phi(1/k))D = \phi(k)D$.
Thus, similarly,
\be\label{eqlamd}
\lambda_1^\uparrow(\Phi_D(G)D^{-1})\ge \phi(\lambda_1^\uparrow(G D^{-1})).
\ee

To combine (\ref{eqlam1}) and (\ref{eqlamd}) into an expression relating the metric distance
$\delta_\infty(\Phi_D(G),D)$ to $\delta_\infty(G,D)$,
we introduce the function
$$
h(x) = \log \phi(\exp(x)).
$$
From $\phi(1/x) = 1/\phi(x)$, we see that $h$ is odd, $h(-x) = -h(x)$.
Moreover, $h$ is monotonously increasing.
Finally, we note that for $-1< p < 1$, $h(x)/x$ achieves its maximum in $x=0$, and
\be\label{h}
\lim_{x\rightarrow 0}\frac{h(x)}{x} = \frac{p}{2^{p+1}-2} =: \beta_p.
\ee

Taking the logarithm of (\ref{eqlam1}) and (\ref{eqlamd}) gives
\beas
y_1:=\log\lambda_1^\downarrow(\Phi_D(G)D^{-1}) &\le& h(\log\lambda_1^\downarrow(G D^{-1})) =: h(x_1), \\
y_2:=\log\lambda_1^\uparrow(\Phi_D(G)D^{-1}) &\ge& h(\log\lambda_1^\uparrow(G D^{-1})) =: h(x_2),
\eeas
where we also introduced some shorthand.
These two inequalities can be combined as $h(x_2)\le y_2\le y_1\le h(x_1)$, showing
that the interval $[y_2,y_1]$ is completely contained in $[h(x_2),h(x_1)]$.
Therefore,
$$
\max(|y_1|,|y_2|) \le \max(|h(x_1)|,|h(x_2)|).
$$
Since $h$ is odd, $|h(x)| = h(|x|)$, and because $h$ is monotonously increasing,
$$
\max(|y_1|,|y_2|) \le \max(h(|x_1|),h(|x_2|)) = h(\max(|x_1|,|x_2|)).
$$
Now the left-hand side is nothing but $\delta_\infty(\Phi_D(G),D)$,
and the right-hand side is $h(\delta_\infty(G,D))$. By (\ref{h})
it finally follows that
$$
\delta_\infty(\Phi_D(G),D) \le \beta_p \delta_\infty(G,D),
$$
which proves that $D=C$ is the only stationary point of $f(D)$.
This finishes the proof of Proposition \ref{prop1}.
%%%%%%%%%%%%%%%%%%%%%%%%%%%%%%%%%%%%%%%%%%%%%%%%%%%%%%%%%%%%%%%%%%%%%%%%%%
\section{Value of $f(D)$ for non-invertible and/or unbounded $D$\label{sec7}}
In this Section we study the behaviour of $f(D)$ for $D$ on the boundary of the PSD cone, that is,
for non-invertible and/or unbounded $D$.
This will result in a proof of Proposition \ref{prop2}.
As mentioned above, this Proposition is used to inductively prove the statement $f(D)\le f(C)$, and relies
on the induction hypothesis that $f(D)\le f(C)$ holds for matrices of lesser dimension.

We consider blocks $C$ and $D$ of size $d\times d$.
Let $P$ be a projector on a $d'$-dimensional subspace of the full $d$-dimensional space,
and let $P^\perp=\id-P$ be the projector on the complementary subspace.

We consider $D$ of the form $D=D'+\epsilon P$,
where $D'$ is bounded and invertible on the complementary subspace ($P^\perp$) and 0 elsewhere.
We study non-invertible $D$ by taking $P$ the projector on the kernel of $D$ and letting $\epsilon$ tend to zero.
Likewise, we study unbounded $D$ by taking $P^\perp$ the projector on the subspace on which $D$ is bounded
and letting $\epsilon$ tend to infinity.

Thus $D^{-1}=D'^{-1}\oplus P/\epsilon$.
Denote $Q:=P^\perp C^2 P^\perp$, $R:=P C^2 P$ and $G'=D'^{-1/2}Q D'^{-1/2}$, thus $G=G'\oplus R/\epsilon$.
Then
\bea
f(D) &=& (\trace(G'+D')^q-\trace G'^q-\trace D'^q) \nonumber\\
&& + (\trace(R/\epsilon+\epsilon P)^q-\trace(R/\epsilon)^q-\trace(\epsilon P)^q). \label{f2}
\eea

We assume validity of the induction hypothesis on the complementary subspace, namely that
$$
\trace(G'+D')^q-\trace G'^q-\trace D'^q
$$
is maximal for $D'=G'$.
Noting that the role of block $C$ in the definition of $f(D)$ is taken up here by $Q^{1/2}$, $D'=G'$ corresponds to $D'=Q^{1/2}$.

We now show that when $q<2$, the second term tends to 0 if $\epsilon$ tends to 0.
By the Lieb-Thirring inequality, and restricting $R^{-1}$ to the subspace of $P$,
$$
\trace(R+\epsilon^2 P)^q = \trace(R(P+\epsilon^2 R^{-1}))^q \le \trace(R^q(P+\epsilon^2 R^{-1})^q).
$$
Since the non-zero eigenvalues of $P+\epsilon^2 R^{-1}$ are all $\ge 1$, we have
$(P+\epsilon^2 R^{-1})^q\le(P+\epsilon^2 R^{-1})^2$, for $q\le 2$, so that also
$$
\trace(R+\epsilon^2 P)^q \le \trace(R^q(P+\epsilon^2 R^{-1})^2).
$$
Hence
\beas
\lefteqn{\trace(R/\epsilon+\epsilon P)^q-\trace(R/\epsilon)^q-\trace(\epsilon P)^q} \\
&=& \epsilon^{-q}(\trace(R+\epsilon^2 P)^q-\trace R^q-\trace(\epsilon^2 P)^q) \\
&\le& \epsilon^{-q}(\trace(R+\epsilon^2 P)^q-\trace R^q) \\
&\le& \epsilon^{-q}(\trace(R^q(P+\epsilon^2 R^{-1})^2)-\trace R^q) \\
&=& \epsilon^{-q}(2\epsilon^2\trace R^{q-1}+\epsilon^4\trace R^{q-2}) \\
&=& 2\epsilon^{2-q}\trace R^{q-1}+\epsilon^{4-q}\trace R^{q-2}.
\eeas
It is easily seen that for values of $q<2$, this tends to 0 if $\epsilon$ does.

The proof that $\trace(R/\epsilon+\epsilon P)^q-\trace(R/\epsilon)^q-\trace(\epsilon P)^q$ tends to 0
if $\epsilon$ tends to infinity is completely similar.

By the induction hypothesis, the first term in (\ref{f2}) obeys the inequality
$$
\trace(G'+D')^q-\trace G'^q-\trace D'^q \le (2^q-2)\trace(Q^{1/2})^q.
$$
Now $Q=P^\perp C^2 P^\perp$ means that in some basis $Q$ is a principal submatrix of $C^2$.
Hence, by eigenvalue interlacing, and by the non-negativity of $C$ and $Q^2$,
$\trace Q^{q/2}\le \trace (C^2)^{q/2} = \trace C^q$, so that
$$
\trace(G'+D')^q-\trace G'^q-\trace D'^q \le f(C).
$$

Combining the two terms proves $f(D)\le f(C)$ for non-invertible/unbounded $D$ with a $d'$-dimensional bounded invertible part,
based on the induction hypothesis $f(D)\le f(C)$ for dimension $d'$.
This finishes the proof of Proposition \ref{prop2}.
%%%%%%%%%%%%%%%%%%%%%%%%%%%%%%%%%%%%%%%%%%%%%%%%%%%%%%%%%%%%%%%%%%%%%%%%%%
\section{Final Remark}
The method used to prove that $G=D$ is the unique solution of (\ref{eq}) can be employed for other matrix equations.
Here we illustrate this for the equation
\be\label{eqfur}
A X^q A = X A^q X,\quad A\ge 0
\ee
and show that $X=A$ is its unique PSD solution when $0\le q<2$.
Again we can use Lemma \ref{lem1} to solve the right-hand side for $X$, giving the equation
$$
X = (A X^q A) \geo A^{-q} = A(X^q \geo A^{-q-2})A.
$$
This defines the map $\Psi_A$:
$$
X\mapsto \Psi_A(X) := A(X^q \geo A^{-q-2})A.
$$
We show that
\be\label{contr2}
\delta_\infty(\Psi_A(X),A) \le (q/2)\delta_\infty(X,A).
\ee
To do so, we consider the log-majorisation version (\cite{ando_hiai}, Theorem 3.1) of Furuta's inequality \cite{furuta}.
Let $\#_\alpha$ denote the $\alpha$-power mean, then
for $A,B\ge 0$, $0<\alpha\le 1$, $p\ge 0$ and $r\le\min(\alpha,\alpha p)$
$$
A^{(1-\alpha)/2} B^\alpha A^{(1-\alpha)/2} \succ_{\log} \left(A^{p-r} \,\,\#_\alpha\,\,
(A^{(1-\alpha)r/2\alpha} B^p A^{(1-\alpha)r/2\alpha}) \right)^{1/p}.
$$
Substituting $A$ by $A^2$, $B$ by $X^{-2}$, $\alpha$ by $1/2$, $p$ by $q/2$, and $r$ by $-1/2$ yields
\beas
A^{1/2} X^{-1} A^{1/2} &\succ_{\log}& (A^{(1+q)/2}(A^{-1-q/2} (X^{-1})^q A^{-1-q/2})^{1/2} A^{(1+q)/2})^{1/(q/2)} \\
&=& (A^{(1+q)/2}(A^{1+q/2} X^q A^{1+q/2})^{-1/2} A^{(1+q)/2})^{1/(q/2)}.
\eeas
From this log-majorisation relation follows directly
that
\beas
\lefteqn{|||\log(A^{1/2} X^{-1} A^{1/2})|||} \\
& \ge & (1/(q/2)) |||\log(A^{(1+q)/2}(A^{1+q/2} X^q A^{1+q/2})^{-1/2} A^{(1+q)/2})|||,
\eeas
for any unitarily invariant norm, hence (\ref{contr2}) indeed holds.
%-----------------------------------------------------------------------
\begin{ack}
This work was supported by The Leverhulme Trust (grant F/07 058/U),
and is part of the QIP-IRC (www.qipirc.org) supported by EPSRC (GR/S82176/0).
The author thanks J.~Eisert for his constructive comments.
\end{ack}
%------------------------------------------------------------- BIBLIOGRAPHY

%%%%%%%%%%%%%%%%%%%%%%%%%%%%%%%%%%%%%%%%%%%%%%%%%%%%%%%%%%%%%%%%%%%
\end{document}